\documentclass[12pt]{article}

\usepackage{amsfonts}

\newtheorem{thm}{Theorem}[section]
\newtheorem{lemma}[thm]{Lemma}
\newtheorem{cor}[thm]{Corollary}
\newtheorem{definition}[thm]{Definition}

\newtheorem{conjecture}{Conjecture}

\newenvironment{remark}{\par\medskip\noindent{\bf Remark.\ }}{\par\smallskip}

\newcommand{\be}{\begin{equation}}
\newcommand{\ee}{\end{equation}}
\newcommand{\openbox}{\leavevmode
  \hbox to8pt{\hfil\vrule\vbox to6pt{\hrule width6pt\vfil\hrule}\vrule}}

\newcommand{\qed}{\hbox to5pt{ } \hfill \openbox\bigskip\medskip}

\newcommand{\Zp}{\mathbb Z _p}

\newcommand{\cF}{\mbox{$\cal F$}}
\newcommand{\cG}{\mbox{$\cal G$}}

\newcommand{\cP}{\mbox{$\cal P$}}

\newcommand{\cQ}{\mbox{$\cal Q$}}

\newcommand{\R}{\mathbb R}

\newcommand{\F}{\mathbb F}

\newcommand{\ve}[1]{\mathbf{#1}}
\title{About F\"uredi's conjecture}
\author{G\'abor Heged\"{u}s
\\{\normalsize  \'Obuda University}
\\{\normalsize B\'ecsi \'ut 96/B, Budapest, Hungary, H-1032}
\\{\normalsize hegedus.gabor@uni-obuda.hu}
}
\begin{document}
\maketitle

\begin{abstract}
Let $t$ be  a non-negative integer and $\mbox{$\cal P$}=\{(A_i,B_i)\}_{1\leq i\leq m}$ be a set-pair family satisfying  $|A_i \cap B_i|\leq t$ for  $1\leq i \leq m$. $\mbox{$\cal P$}$ is called {\em strong Bollob\'as $t$-system}, if $|A_i\cap B_j|>t$ for all $1\leq i\neq  j \leq m$. 

F\"uredi conjectured  the following nice generalization of Bollob\'as' Theorem:\\
Let $t$ be  a non-negative integer. Let $\mbox{$\cal P$}=\{(A_i,B_i)\}_{1\leq i\leq m}$ be a strong Bollob\'as  $t$-system. 
Then
$$
\sum_{i=1}^m \frac{1}{{|A_i|+|B_i|-2t \choose |A_i|-t}}\leq 1.
$$
We confirmed  the following special case of F\"uredi's conjecture along with some more results of similar flavor.\\
Let $t$ be  a non-negative integer. Let   $\mbox{$\cal P$}=\{(A_i,B_i)\}_{1\leq i\leq m}$ denote a strong Bollob\'as  $t$-system. Define $a_i:=|A_i|$ and $b_i:=|B_i|$ for each $i$. 
Assume that there exists a positive integer $N$ such that $a_i+b_i=N$ for each $i$. Then 
$$
\sum_{i=1}^m \frac{1}{{a_i+b_i-2t \choose a_i-t}}\leq 1.
$$
\end{abstract}
\medskip


\section{Introduction}

First we introduce some notation and  give some basic definitions. 
Let $n$ be a positive integer and $[n]$ stand for the set $\{1,2,
\ldots, n\}$. For an integer $0\leq d\leq n$ we denote by
${[n] \choose d}$ the family of all  $d$ element subsets of $[n]$.

\begin{definition}
Let $\cP=\{(A_i,B_i)\}_{1\leq i\leq m}$, where $A_i,B_i\subseteq [n]$, $A_i\cap B_i=\emptyset$ for all $1\leq i\leq m$. Then $\cP$ is called a {\em Bollob\'as  system} (or strong Bollob\'as  system) if  $A_i\cap B_j\ne \emptyset$ holds whenever $i\ne j$. Also, $\cP$ is called a {\em skew Bollob\'as  system} if  $A_i\cap B_j\ne \emptyset$ is only required for all $1\leq i<j \leq m$.
\end{definition}

\begin{definition}
Let $t$ be  a non-negative integer and $\cP=\{(A_i,B_i)\}_{1\leq i\leq m}$ be a set-pair family satisfying  $|A_i \cap B_i|\leq t$ for  $1\leq i \leq m$. Then $\cP$ is called {\em Bollob\'as $t$-system} ({\em skew Bollob\'as $t$-system}), if $|A_i\cap B_j|>t$ for all $1\leq i\neq  j \leq m$ ($1\leq i<  j \leq m$), respectively. 
\end{definition}

\begin{definition}
Let $\F$ be a fixed field. Let $\cP=\{(U_i,V_i)\}_{1\leq i\leq m}$, where $U_i,V_i\leq W$ are subspaces of an $n$-dimensional  vector space $W$ over $\F$. Suppose that  $\dim(U_i\cap V_i)=0$ for all $1\leq i\leq m$. Then $\cP$ is called a {\em Bollob\'as subspace system} (or strong Bollob\'as subspace  system) if  $\dim(U_i\cap V_j)\neq 0$ holds whenever $i\ne j$. Also, $\cP$ is called a {\em skew Bollob\'as subspace system} if  $\dim(U_i\cap V_j)\neq 0$  is only required for all $1\leq i<j \leq m$.
\end{definition}

\begin{definition}
Let $t$ be  a non-negative integer. Let $\F$ be a fixed field. Let $\cP=\{(U_i,V_i)\}_{1\leq i\leq m}$, where $U_i,V_i\leq W$ are subspaces of an $n$-dimensional  vector space $W$ over $\F$.  Suppose that  $\dim(U_i\cap V_i)\leq t$ for all $1\leq i\leq m$. Then $\cP$ is called a {\em Bollob\'as subspace $t$-system} (or strong Bollob\'as subspace $t$-system) if  $\dim(U_i\cap V_j) >t$ holds whenever $i\ne j$. Also, $\cP$ is called a {\em skew Bollob\'as subspace $t$-system} if  $\dim(U_i\cap V_j)> t$  is only required for all $1\leq i<j \leq m$.
\end{definition}

Bollob\'as proved  the following remarkable result in extremal combinatorics in \cite{B}.

\begin{thm} (\cite{B}) \label{Bollobas}
Let  $\cP=\{(A_i,B_i)\}_{1\leq i\leq m}$ be a strong Bollob\'as  system. Then
\begin{equation}  \label{Boll}
\sum_{i=1}^m \frac{1}{{|A_i|+|B_i| \choose |A_i|}}\leq 1.
\end{equation}
\end{thm}

It is straightforward to verify that $\cF=\{F_1, \ldots ,F_m\}$ is an antichain iff 
$\cP=\{(F_i,[n]\setminus F_i)\}_{1\leq i\leq m}$ is  a strong Bollob\'as system. Thus Theorem (\ref{Boll}) implies 
\\
{\em LYM-inequality} (see \cite{Lu}, \cite{M}, \cite{Y}). Suppose that $\cF=\{F_1, \ldots ,F_m\}$ is an antichain. Then 
\begin{equation}  \label{LYM}
\sum_{i=1}^m  \frac{1}{{n\choose |F_i|}}\leq 1.
\end{equation}

Lov\'asz gave the following skew uniform version of Bollob\'as' Theorem for subspaces in \cite{L1}. He introduced  a new proof method, the tensor product method to prove this version.

\begin{thm} \label{Lovasz}
Let $\F$ be an arbitrary field.  Let  $U_1, \ldots ,U_m$ be $r$-dimensional and $V_1, \ldots ,V_m$ be $s$-dimensional subspaces of a linear  space $W$ over the field $\F$. Assume that 
  $\cP=\{(U_i,V_i)\}_{1\leq i\leq m}$ is a skew Bollob\'as subspace system.
Then
$$
m\leq {r+s \choose r}.
$$
\end{thm}

F\"uredi generalized Lov\'asz' results in \cite{F} and he obtained a nice threshold variant of Theorem  \ref{Lovasz}. 

\begin{thm} \label{Furedi1}
Let $t$ be  a non-negative integer. Let $\F$ be an arbitrary field.  Let  $U_1, \ldots ,U_m$ be $r$-dimensional and $V_1, \ldots ,V_m$ be $s$-dimensional subspaces of a linear  space $W$ over the field $\F$.  Assume that  $\cP=\{(U_i,V_i)\}_{1\leq i\leq m}$ is a skew Bollob\'as subspace $t$-system.
Then
$$
m\leq {r+s-2t \choose r-t}.
$$
\end{thm}

F\"uredi derived the following combinatorial result from Theorem 
\ref{Furedi1}.

\begin{cor} \label{Furedi2}
Let $t$ be  a non-negative integer. Let $A_1, \ldots ,A_m$ be $r$-element sets and  $B_1, \ldots ,B_m$ be $s$-element sets. Assume that  $\cP=\{(A_i,B_i)\}_{1\leq i\leq m}$ is a skew Bollob\'as  $t$-system. 
Then
$$
m\leq {r+s-2t \choose r-t}.
$$
\end{cor}

F\"uredi conjectured  the following nice generalization of Bollob\'as' Theorem  in \cite{F}.
\begin{conjecture} \label{Fconj}
Let $t$ be  a non-negative integer. Let $\cP=\{(A_i,B_i)\}_{1\leq i\leq m}$ be a strong Bollob\'as  $t$-system. 
Then
\begin{equation}  \label{Fur}
\sum_{i=1}^m \frac{1}{{|A_i|+|B_i|-2t \choose |A_i|-t}}\leq 1.
\end{equation}
\end{conjecture}
\begin{remark}
It is easy to verify that Conjecture \ref{Fconj} is sharp, if it is true.
Namely let $a,b,t>0$ positive integers and let $n:=a+b+t$. Let $T=\{a+b+1, \ldots ,n\}$, then $|T|=t$. Consider the set system $\cG={[a+b] \choose a}$. 

Consider the system $\cP:=\{(A\cup T, ([a+b]\setminus A) \cup T)\}_{A\in \cG}$. Clearly $m=|\cG|={a+b \choose a}$.  We can check easily that $\cP$ is a strong Bollob\'as  $t$-system and 
$$
\sum_{i=1}^m \frac{1}{{a+b \choose a}}= 1.
$$
\end{remark}

Zhu made some progress in the solution of  Conjecture \ref{Fconj} and  he proved the following weaker result in \cite{Z}.

\begin{thm} \label{Zhu}
Let $A_1, \ldots ,A_m$ and $B_1, \ldots ,B_m$ be finite sets such that $|A_i \cap B_i|=t$ for  each $1\leq i \leq m$ and $|A_i \cap B_j|>t$ for each $i\neq j$. Then
$$
\sum_{i=1}^m \frac{1}{{|A_i|+|B_i|-t \choose |B_i|-t}}\leq 1.
$$
\end{thm}
 

One of our main result gives a natural inequality for  non-uniform skew Bollob\'as subspace $t$-systems.

\begin{thm} \label{main1}
Let $t$ be  a non-negative integer. Let   $\cP=\{(U_i,V_i)\}_{1\leq i\leq m}$ denote a skew Bollob\'as subspace $t$-system over the real field $\R$. Define $u_i:=\dim(U_i)$ and $v_i:=\dim(V_i)$ for each $i$. 
Assume that $u_1\leq u_2\leq  \ldots 	\leq u_m$ and $v_1\geq v_2\geq  \ldots 	\geq v_m$.
Then
\begin{equation}  \label{main_eq}
\sum_{i=1}^m \frac{1}{{u_i+v_i-2t \choose u_i-t}}\leq 1.
\end{equation}
\end{thm}

\begin{remark}
Clearly Theorem \ref{Furedi1} is a special case of  Theorem \ref{main1} in the real case $\F=\R$.
\end{remark}

We can derive easily the following consequences of Theorem \ref{main1}.
 
\begin{cor} \label{main2}
Let $t$ be  a non-negative integer. Let   $\cP=\{(U_i,V_i)\}_{1\leq i\leq m}$ denote a strong Bollob\'as subspace $t$-system over the real field $\R$. Define $u_i:=\dim(U_i)$ and $v_i:=\dim(V_i)$ for each $i$. 
Assume that there exists a  positive integer $N$ such that $u_i+v_i=N$ for each $i$. Then
\begin{equation}  \label{main_eq2}
\sum_{i=1}^m \frac{1}{{u_i+v_i-2t \choose u_i-t}}\leq 1.
\end{equation}
\end{cor}

\begin{cor} \label{main3}
Let $t$ be  a non-negative integer. Let   $\cP=\{(U_i,V_i)\}_{1\leq i\leq m}$ denote a strong Bollob\'as subspace $t$-system over the real field $\R$. Define $u_i:=\dim(U_i)$  for each $i$  and let $v:=\max_i\dim(V_i)$. 
Then 
\begin{equation}  \label{main_eq3}
\sum_{i=1}^m \frac{1}{{u_i+v-2t \choose u_i-t}}\leq 1.
\end{equation}
\end{cor}

We derive also the following special cases of 
F\"uredi's conjecture.

\begin{cor} \label{main4}
Let $t$ be  a non-negative integer. Let   $\cP=\{(A_i,B_i)\}_{1\leq i\leq m}$ denote a skew Bollob\'as  $t$-system. Define $a_i:=|A_i|$ and $b_i:=|B_i|$ for each $i$. 
Assume that $a_1\leq a_2\leq  \ldots 	\leq a_m$ and $b_1\geq b_2\geq  \ldots 	\geq b_m$.
Then
\begin{equation}  \label{main_eq4}
\sum_{i=1}^m \frac{1}{{a_i+b_i-2t \choose a_i-t}}\leq 1.
\end{equation}
\end{cor}

\begin{remark}
If $t=0$, then Corollary \ref{main4} is not true without the assumption $a_1\leq a_2\leq  \ldots 	\leq a_m$ and $b_1\geq b_2\geq  \ldots 	\geq b_m$ (see e.g. the counterexample  \cite{FH} Example 1).
\end{remark}

\begin{cor} \label{main5}
Let $t$ be  a non-negative integer. Let   $\cP=\{(A_i,B_i)\}_{1\leq i\leq m}$ denote a strong Bollob\'as  $t$-system. Define $a_i:=|A_i|$ and $b_i:=|B_i|$ for each $i$. 
Assume that there exists a  positive integer $N$  such that $a_i+b_i=N$ for each $i$. Then
\begin{equation}  \label{main_eq5}
\sum_{i=1}^m \frac{1}{{a_i+b_i-2t \choose a_i-t}}\leq 1.
\end{equation}
\end{cor}

\begin{cor} \label{main6}
Let $t$ be  a non-negative integer. Let   $\cP=\{(U_i,V_i)\}_{1\leq i\leq m}$ denote a strong Bollob\'as  $t$-system. Define $a_i:=|A_i|$  for each $i$  and let $b:=\max_i\ |B_i|$. 
Then 
\begin{equation}  \label{main_eq6}
\sum_{i=1}^m \frac{1}{{a_i+b-2t \choose a_i-t}}\leq 1.
\end{equation}
\end{cor}

The proof of Theorem \ref{main1} combines  tensor product methods and F\"uredi's method, which he applied to prove Theorem \ref{Furedi1} in \cite{F}. 

We present our proofs in Section 2. In Section 3 we raise some natural conjectures which imply  F\"uredi's conjecture.

\section{Proofs}

Our proof is completely based on Scott and Wilmer's following result (see \cite{SW} Theorem 4.5).

\begin{thm} \label{Scott}
Let   $\cP=\{(U_i,V_i)\}_{1\leq i\leq m}$ denote a skew Bollob\'as subspace system over the real field $\R$. Define $u_i:=\dim(U_i)$ and $v_i:=\dim(V_i)$ for each $i$. 
Assume that $u_1\leq u_2\leq  \ldots 	\leq u_m$ and $v_1\geq v_2\geq  \ldots 	\geq v_m$.
Then
\begin{equation}  \label{scott}
\sum_{i=1}^m \frac{1}{{u_i+v_i \choose u_i}}\leq 1.
\end{equation}
\end{thm}

We use in our proof the following well-known Lemma  (see e.g.   \cite{FT} Lemma 26.14). 
\begin{lemma} \label{subs_gen_pos}
Let $n\geq k$ and $t=n-k$. Let $\F$ be an arbitrary infinite field. Let $V$ denote an $n$-dimensional vector space. Let $W_1, \ldots ,W_m$ be subspaces of $V$ with $\dim(W_i)<n$ for each $i$. Then there exists a $k$-dimensional subspace $V'$ such that 
$$
\dim(W_i\cap V')=\mbox{\rm max}(\dim(W_i)-t,0)
$$ 
for each $1\leq i\leq m$.          
\end{lemma}

\begin{remark}
This subspace $V'$ guaranteed by Lemma  \ref{subs_gen_pos} is called to be in {\em general position} with respect to the subspaces $W_1, \ldots ,W_m$.
\end{remark}

{\bf Proof of Theorem \ref{main1}:}

It follows from Lemma \ref{subs_gen_pos} that there exists a subspace of co-dimension $t$ in general position with respect to the subspaces $U_i$, $V_i$ and $U_i\cap V_i$ (here $1\leq i\leq m$). Let $W_0$ denote this subspace. Clearly $\dim(U_i\cap W_0)=\dim(U_i)-t$ and $\dim(V_i\cap W_0)=\dim(V_i)-t$ for each $i$. 

Consequently
\begin{itemize}
\item[(i)] $\dim(U_i \cap V_i\cap W_0)= 0$ for each $1\leq i \leq m$;
\item[(ii)] $\dim(U_i\cap V_j\cap W_0)>0$ whenever $i< j$ ($1\leq i, j \leq m$);
\item[(iii)] $\dim(U_i\cap W_0)\leq \dim(U_j\cap W_0)$ and $\dim(V_i\cap W_0)\geq \dim(V_j\cap W_0)$ whenever $i< j$ ($1\leq i, j \leq m$).
\end{itemize}

Hence the conditions (i), (ii) and (iii) guarantee that for the subspaces $U_i\cap W_0$, $V_i\cap W_0$ and $W_0$ the conditions of Theorem  \ref{Scott} hold. Therefore  we can apply Theorem \ref{Scott} for the subspaces $U_i':=U_i\cap W_0$, $V_i':=V_i\cap W_0$ and $W_0$, consequently
$$
\sum_{i=1}^m \frac{1}{{u_i+v_i-2t \choose u_i-t}}=\sum_{i=1}^m \frac{1}{{u'_i+v'_i \choose u'_i}}\leq 1,
$$
where $u'_i:=\dim(U'_i)=\dim(U_i)-t$ and $v'_i:=\dim(V'_i)=\dim(V_i)-t$ for each $i$. \qed

{\bf Proof of Corollary \ref{main2}:}

Let $\cP=\{(U_i,V_i)\}_{1\leq i\leq m}$ be  a {\em strong} Bollob\'as subspace $t$-system over the real field. Then there exists a permutation $\pi\in S_m$ such that $u_{\pi(i)}\leq u_{\pi(j)}$ for each $i<j$, where $u_i:= \dim(U_i)$ and $v_i:= \dim(V_i)$ for each $i$. But $u_i+v_i=N$ by assumption, which implies that  $v_{\pi(i)}\geq v_{\pi(j)}$ for each $i<j$. Define $U'_i:=U_{\pi(i)}$,  $V'_i:=V_{\pi(i)}$, $u'_i:=u_{\pi(i)}$ and $v'_i:=v_{\pi(i)}$ for each $i$. Then we can apply 
Theorem \ref{main1} for the {\em skew} Bollob\'as subspace $t$-system 
$\cP':=\{(U'_i,V'_i)\}_{1\leq i\leq m}$ and we get that 
$$
\sum_{i=1}^m \frac{1}{{u_i+v_i-2t \choose u_i-t}}=\sum_{i=1}^m \frac{1}{{u'_i+v'_i-2t \choose u'_i-t}}\leq 1,
$$
which was to be proved. \qed

{\bf Proof of Corollary \ref{main3}:}

Let   $\cP=\{(U_i,V_i)\}_{1\leq i\leq m}$ be a {\em strong} Bollob\'as subspace $t$-system over the real field $\R$. Define $u_i:=\dim(U_i)$  for each $i$  and $v:=\max_i\dim(V_i)$. 

Let $1\leq i\leq m$ be a fixed index. It is easy to verify that there exists a $W_i$ subspace such that $V_i\leq W_i$, $\dim(W_i)=v$ and $U_i\cap V_i=U_i\cap W_i$. 

Consider the new subspace system  $\cQ:=\{(U_i,W_i)\}_{1\leq i\leq m}$. Then $\dim(U_i\cap W_i)=\dim(U_i\cap V_i)\leq t$ for each $i$ and $\dim(U_i\cap W_j)\geq \dim(U_i\cap V_j)> t$ for each $1\leq i\neq j\leq m$, hence $\cQ$ is a strong Bollob\'as subspace $t$-system over the real field. 

Let $w_i:=\dim(W_i)$ for each $i$. Clearly $w_1=w_2=\ldots =w_m=v$. After an appropriate permutation we can suppose without lost of generality that $u_1\leq u_2\leq  \ldots 	\leq u_m$.  Hence we can apply Theorem \ref{main1} to the strong Bollob\'as subspace $t$-system $\cQ$ and we get that
$$
\sum_{i=1}^m \frac{1}{{u_i+v-2t \choose u_i-t}}\leq 1.
$$
\qed

\section{Concluding remarks}

We can raise  the following natural  conjectures.

\begin{conjecture} \label{Hconj1}
Let   $\cP=\{(U_i,V_i)\}_{1\leq i\leq m}$ denote a strong Bollob\'as subspace system over the real field $\R$. Define $u_i:=\dim(U_i)$    and $v_i:=\dim(V_i)$ for each $i$. 
Then 
\begin{equation}  \label{main_eq9}
\sum_{i=1}^m \frac{1}{{u_i+v_i \choose u_i}}\leq 1.
\end{equation}
\end{conjecture}
\begin{remark}
Let $a,b>0$ be positive integers and define $n:=a+b$. Denote by $\{\ve e_1, \ldots , \ve e_n\}$ the standard basis of $\R^n$. For each $I\subseteq [n]$ define the subspace $V(I)\leq \R^n$ as the subspace generated by the set $\{ \ve e_i:~ i\in I\}$. 

Consider the subspace system $\cP:=\{(V(I), V([n]\setminus I) )\}_{I\in {[a+b]\choose a}}$. Then it is easy to verify
 that this system $\cP$ shows that if Conjecture \ref{Hconj1} is true, then it is sharp.
\end{remark}

\begin{conjecture} \label{Hconj2}
Let $t$ be  a non-negative integer. Let   $\cP=\{(U_i,V_i)\}_{1\leq i\leq m}$ denote a strong Bollob\'as subspace $t$-system over the real field $\R$. Define $u_i:=\dim(U_i)$    and $v_i:=\dim(V_i)$ for each $i$. 
Then 
\begin{equation}  \label{main_eq10}
\sum_{i=1}^m \frac{1}{{u_i+v_i-2t \choose u_i-t}}\leq 1.
\end{equation}
\end{conjecture}

Finally we note that we  can easily derive  the following result using F\"uredi's method.
\begin{thm} Conjecture \ref{Hconj1} is equivalent with   Conjecture \ref{Hconj2}.
\end{thm}
\begin{remark}
It is easy to verify that Conjecture \ref{Hconj2} implies Conjecture  \ref{Fconj}.
\end{remark}


\end{document}